\documentclass[12pt]{article}
\usepackage{amssymb}
\usepackage{amsfonts}
\usepackage{amsmath}
\usepackage{amsbsy}
\newcommand{\be}{\begin{equation}}
\newcommand{\ee}{\end{equation}}
\newcommand{\bea}{\begin{eqnarray}}
\newcommand{\eea}{\end{eqnarray}}
\newcommand{\ba}{\begin{array}}
\newcommand{\ea}{\end{array}}
\newcommand{\R}{I\!\!R}

\newcommand{\bc}{\begin{center}}
\newcommand{\ec}{\end{center}}
\newcommand{\ben}{\begin{enumerate}}
\newcommand{\een}{\end{enumerate}}
\newcommand{\bfi}{\begin{figure}}
\newcommand{\efi}{\end{figure}}

\newcommand{\bq}{\begin{quote}}
\newcommand{\eq}{\end{quote}}
\newcommand{\bqu}{\begin{quotation}}
\newcommand{\equ}{\end{quotation}}
\newenvironment{emphit}{\begin{itemize}}{\end{itemize}}
\newcommand{\bemp}{\begin{emphit}}
\newcommand{\eemp}{\end{emphit}}

\newcommand{\bt}{\begin{tabular}}
\newcommand{\et}{\end{tabular}}

\newtheorem{myth}{Theorem}[section]
\newtheorem{mylem}{Lemma}[section]
\newtheorem{mycor}{Corollary}[section]

\newtheorem{mydef}{Definition}[section]
\newtheorem{myrem}{Remark}[section]

\begin{document}
\date{}
\title{On the factor set of code loops
}
\author{T\`em\'it\d{\'o}p\d{\'e} Gb\d{\'o}l\'ah\`an Ja\'iy\'e\d ol\'a
\thanks{All correspondence to be addressed to this author.} \\
Department of Mathematics,\\
Obafemi Awolowo University, Ile Ife, Nigeria.\\
jaiyeolatemitope@yahoo.com,~tjayeola@oauife.edu.ng} \maketitle

\begin{abstract}
A Code loop on a binary linear code that is doubly even with a
factor set is shown to be a central loop, conjugacy closed loop,
Burn loop and extra loop. General forms of the identities that
define the factor set of a code are deduced.
\end{abstract}

\section{Introduction}
Let $L$ be a loop and $Z\le Z(L)$. If $L/Z\cong C$ where $C$ is an
abelian group then, $L\cong L(\phi )$ where $\phi~:~C\times C\to F$,
$F$ been a field and $L(\phi )=F\times C$ with binary operation
$\ast$~ such that for all $(a,u),(b,v)\in L(\phi)$ ;
\begin{equation}\label{eq:1}
(a,u)\ast (b,v)=(a+b+\phi (u,v),u+v)
\end{equation}
and $(L(\phi ),\ast)$ is a loop. The fact that $(L(\phi ),\ast)$ is
a loop is shown in \cite{phd21}.

Let $V$ be a finite dimensional vector space over the field $F$ and
$C\le V$. With $|F|=2$ such that $F=\mathbb{Z}_2=\{0,1\}$, $C$ is
called a binary linear code and its elements (vectors) are called
code words. For all $u,v\in C$, let $||u||$ denote the number of
non-zero co-ordinates in $u$ (this is simply equal to the 'norm
squared' of $u$ if $V={\R}^{n}$), let $u\cdot v$ denote the number
of the corresponding co-ordinates of $u$ and $v$ that are non-zero
(this is simply equal to the 'inner product' of $u$ and $v$ if
$V={\R}^n$) and let $c(u,v,w)$ represent the number of corresponding
co-ordinates of $u$, $v$ and $w$ that are non-zero.

$C$ is called a double even code if :
\begin{enumerate}
\item $u\cdot v$ is even
\item $||u||$ is divisible by $4$.
\end{enumerate}
$\phi$ is called a factor set on $C$ if it satisfy the conditions :
\begin{equation}\label{eq:2}
\phi (u,u)\equiv \frac{||u||}{4}\bmod{2}
\end{equation}
\begin{equation}\label{eq:3}
\phi (u,v)+\phi (v,u)\equiv \frac{u\cdot v}{2}\bmod{2}
\end{equation}
\begin{equation}\label{eq:4}
\phi (u,v)+\phi (u+v,w)+\phi (v,w)+\phi (u,v+w)\equiv
c(u,v,w)\bmod{2}
\end{equation}
As mentioned in \cite{phd14}, for all $u\in C$, $\phi (u,0)=\phi
(0,u)=0$.

For more on the varieties of loops and their general properties,
readers can check \cite{phd41},\cite{phd39}, \cite{phd49},
\cite{phd42} and \cite{phd3}.

Code loops were first shown to be Moufang loops by Griess
\cite{phd66} and was later shown to be Moufang loops with unique
non-identity square in \cite{phd20}. They have also been of interest
to Chein and Goodaire \cite{phd115} and Vojt\v echovsk\'y
\cite{phd69}, \cite{phd67}, \cite{phd68} and \cite{phd21}. Here we
shall approach them in the style of \cite{phd20}.

A Code loop on a binary linear code that is doubly even with a
factor set is shown to be a central loop, conjugacy closed loop,
Burn loop and extra loop. General forms of the identities that
define the factor set of a code are deduced.

In all, LB, RB, LC, RC, C, LCC, RCC, and CC represent left Bol,
right Bol, left central, right central, central, left conjugacy
closed, right conjugacy closed and conjugacy closed respectively.

\section{Preliminaries}
\begin{mydef}\label{1:0.1}
Let $(G,\circ )$ be a loop. $G$ is a left(right) Burn loop if and
only if it is a LB(RB)-loop and a LCC(RCC)-loop. It is called a Burn
loop if and only if it is both a Moufang loop and a CC-loop.
\end{mydef}

\begin{mydef}\label{1:1}
Let $\phi$ be the factor set of a binary linear code $C$ that is
doubly even. $\phi$ is called a left weakly linear (lwl) factor set
on $C\Leftrightarrow \phi (u+v,u)=\phi (u,u)+\phi
(v,u)~\forall~u,v\in C$.

Similarly, $\phi$ is called a right weakly linear (rwl) factor set
on $C\Leftrightarrow \phi (u,u+v)=\phi (u,u)+\phi
(u,v)~\forall~u,v\in C$.

$\phi$ is called a weakly linear (wl) factor set on
$C\Leftrightarrow \phi$ is both a lwl and rwl factor set on $C$.
\end{mydef}

\begin{mylem}\label{1:2}
Let $\phi$ be the factor set of a binary linear code $C$ that is
doubly even.
\begin{enumerate}
\item If $\phi$ is lwl, then $\phi (u+v,u)\equiv [\phi (u,u)+\frac{u\cdot
v}{2}+\phi(u,v)]\bmod{2}$.
\item If $\phi$ is rwl, then $\phi (u,u+v)\equiv [\phi (u,u)+\frac{u\cdot
v}{2}+\phi(v,u)]\bmod{2}$.
\end{enumerate}
\end{mylem}
{\bf Proof}\\
In a doubly even code $C$ with factor set $\phi$ : $\phi (u,v)+\phi
(v,u)\equiv \frac{u\cdot v}{2}\bmod{2}$. Thus ;
\begin{equation}\label{eq:5}
\phi (u,v)\equiv \Big(\frac{u\cdot v}{2}+\phi (v,u)\Big)\bmod{2}
\end{equation}
\begin{equation}\label{eq:6}
\phi (v,u)\equiv \big(\frac{u\cdot v}{2}+\phi (u,v)\big)\bmod{2}
\end{equation}
\begin{enumerate}
\item Since $\phi$ is lwl, then using
(\ref{eq:6}), $\phi (u+v,u)\equiv [\phi (u,u)+\frac{u\cdot
v}{2}+\phi(u,v)]\bmod{2}$.
\item Since $\phi$ is rwl, then using (\ref{eq:5}),
$\phi (u,u+v)\equiv [\phi (u,u)+\frac{u\cdot
v}{2}+\phi(v,u)]\bmod{2}$.
\end{enumerate}

\begin{mycor}\label{1:2.1}
Let $\phi$ be the factor set of a binary linear code $C$ that is
doubly even. If $\phi$ is wl, then $\phi (u+v,u)+\phi (u,u+v)\equiv
\frac{u\cdot v}{2}\bmod{2}$.
\end{mycor}
{\bf Proof}\\
This follows by Lemma~\ref{1:2}.

\begin{mycor}\label{1:3}
Let $\phi$ be the factor set of a binary linear code $C$ that is
doubly even . The difference between $\phi$ been rwl and lwl is
equal to half of the number of corresponding co-ordinates that are
non-zero in any two code words.
\end{mycor}
{\bf Proof}\\
This follows from the fact in Corollary~\ref{1:2.1} above. $\phi
(u+v,u)-\phi (u,u+v)\equiv \frac{u\cdot v}{2}\bmod{2}$ and $\phi
(u,u+v)-\phi (u+v,u)\equiv \frac{u\cdot v}{2}\bmod{2}$. Hence the
proof.

\begin{mylem}\label{1:4}
Let $\phi$ be the factor set of a binary linear code $C$ that is
doubly even. If $\phi$ is both rwl and lwl, then : $\phi
(u,u+v)=\phi (u+v,v)\Leftrightarrow \phi (u,v)=\phi
(v,u)~\forall~u,v\in C$.
\end{mylem}
{\bf Proof}\\
With $\phi$ been rwl and lwl ; $\phi (u,u+v)=\phi
(u+v,v)\Leftrightarrow \phi(u,u)+\phi (u,v)=\phi (u,u)+\phi
(v,u)\Leftrightarrow \phi (u,v)=\phi (v,u)$.

\begin{mycor}\label{1:4.2}
Let $\phi$ be the factor set of a binary linear code $C$ that is
doubly even. If $\phi$ is wl  then : $\phi (u,u+v)=\phi
(u+v,v)\Leftrightarrow \phi (u,v)=\phi (v,u)~\forall~u,v\in C$.
\end{mycor}
{\bf Proof}\\
This follows from Lemma~\ref{1:4}.
\newpage

{\huge Main Results}

\section{Code Loops and Central Loops}
\begin{myth}\label{1:5}
Let $L(\phi)$ be a code loop on a binary linear code $C$ that is
doubly even with factor set $\phi$. $L(\phi)$ is an LC-loop if and
only if any of the following equivalent conditions is true for all
$u,v,w\in C$.
\begin{enumerate}
\item $A_1(u,v)=\phi (u,u+v)+\phi (u,v)+\phi (u,u)=0$.
\item $A_2(u,v,w)=\phi (u,v)+\phi (u,u+v)+\phi (u,v+w)+\phi
(u,u+v+w)=0$.
\item $A_3(u,v,w)=\phi (u,u)+\phi (u,v+w)+\phi (u,u+v+w)=0$
\end{enumerate}
\end{myth}
{\bf Proof}\\
\begin{enumerate}
\item $L(\phi)$ is a LC-loop $\Leftrightarrow xx\ast yz=(x\ast xy)z~\forall~x=(a,u),y=(b,v),z=(c,w)\in L(\phi)$. $(a,u)(a,u)\ast (b,v)(c,w)=[(a,u)\ast (a,u)(b,v)](c,w)\Leftrightarrow (b+c+\phi
(u,u)+\phi (v,w),v+w)=(b+c+\phi (u,v)+\phi (u,u+v)+\phi
(v,w),v+w)\Leftrightarrow A_1(u,v)=\phi (u,u+v)+\phi (u,v)+\phi
(u,u)=0$.
\item $L(\phi)$ is a LC-loop $\Leftrightarrow (x\ast xy)z=x(
x\ast yz)~\forall~x=(a,u),y=(b,v),z=(c,w)\in L(\phi)$.
$[(a,u)\ast(a,u)(b,v)](c,w)=(a,u)[(a,u)\ast(b,v)(c,w)]\Leftrightarrow
[(a,u)\ast(a+b+\phi (u,v),u+v)](c,w)=(a,u)[(a,u)\ast (b+c+\phi
(v,w),v+w)]\Leftrightarrow (b+c+\phi (u,v)+\phi (u,u+v)+\phi
(v,w),v+w)=(b+c+\phi (v,w)+\phi (u,v+w)+\phi
(u,u+v+w),v+w)\Leftrightarrow A_2(u,v,w)=\phi (u,v)+\phi
(u,u+v)+\phi (u,v+w)+\phi (u,u+v+w)=0$.
\item $L(\phi)$ is a LC-loop $\Leftrightarrow (xx\ast y)z=x(
x\ast yz)~\forall~x=(a,u),y=(b,v),z=(c,w)\in L(\phi)$.
$[(a,u)(a,u)\ast
(b,v)](c,w)=(a,u)[(a,u)\ast(b,v)(c,w)]\Leftrightarrow (b+c+\phi
(u,u)+\phi (v,w),v+w)=(b+c+\phi (v,w)+\phi (u,v+w)+\phi
(u,u+v+w),v+w)\Leftrightarrow A_3(u,v,w)=\phi (u,u)+\phi
(u,v+w)+\phi (u,u+v+w)=0$.
\end{enumerate}

\begin{mycor}\label{1:6}
Let $L(\phi)$ be a code loop on a binary linear code $C$ that is
doubly even with factor set $\phi$. $L(\phi)$ is an LC-loop if and
only if $\phi$ is a rwl factor set.
\end{mycor}
{\bf Proof}\\
Following Theorem~\ref{1:5}, $A_1(u,v)=0\Leftrightarrow \phi
(u,u+v)=\phi (u,u)+\phi (u,v)~\forall~u,v\in C$.

\begin{mycor}\label{1:7}
Let $L(\phi)$ be a code loop on a binary linear code $C$ that is
doubly even with factor set $\phi$. If $L(\phi)$ is an LC-loop ;
\begin{enumerate}
\item $\phi (u,u+v)\equiv [\phi (u,u)+\frac{u\cdot
v}{2}+\phi(v,u)]\bmod{2}$.
\item The difference between $\phi$ been rwl and lwl is $\frac{u\cdot
v}{2}\bmod{2}~\forall~u,v\in C$.
\end{enumerate}
\end{mycor}
{\bf Proof}\\
These follows from Lemma~\ref{1:2} and Corollary~\ref{1:3}.

\begin{myth}\label{1:8}
Let $L(\phi)$ be a code loop on a binary linear code $C$ that is
doubly even with factor set $\phi$. $L(\phi)$ is an RC-loop if and
only if any of the following equivalent conditions is true for all
$u,v,w\in C$.
\begin{enumerate}
\item $B_1(u,w)=\phi (u,u)+\phi (w,u)+\phi (w+u,u)=0$.
\item $B_2(u,w,v)=\phi (w,u)+\phi (w+u,u)+\phi (v+w,u)+\phi
(v+w+u,u)=0$.
\item $B_3(u,v,w)=\phi (u,u)+\phi (v+w,u)+\phi (v+w+u,u)=0$.
\end{enumerate}
\end{myth}
{\bf Proof}\\
\begin{enumerate}
\item $L(\phi)$ is a RC-loop $\Leftrightarrow yz\ast xx=y(zx\ast x)~\forall~x=(a,u),y=(b,v),z=(c,w)\in L(\phi)$. $(b,v)(c,w)\ast(a,u)(a,u))=(b,v)[(c,w)(a,u)\ast(a,u)]\Leftrightarrow (b+c+\phi
(v,w),+\phi (u,u),v+w)=(b+c+\phi (w,u)+\phi (w+u,u)+\phi
(v,w),v+w)\Leftrightarrow B_1(u,w)=\phi (u,u)+\phi (w,u)+\phi
(w+u,u)=0$.
\item $L(\phi)$ is a RC-loop $\Leftrightarrow (yz\ast x)x
=y(zx\ast x)~\forall~x=(a,u),y=(b,v),z=(c,w)\in L(\phi)$.
$[(b,c)(c,w)\ast(a,u)](a,u)=(b,v)[(c,w)(a,u)\ast
(a,u)]\Leftrightarrow (b+c+\phi
(v,w),v+w)(a,u)\ast(a,u)=(b,v)\ast(c+a+\phi
(w,u),w+u)(a,u)\Leftrightarrow (b+c+\phi (v,w)+\phi (v+w,u)+\phi
(v+w+u,u),v+w)=(b+c+\phi (w,u)+\phi (w+u,u)+\phi
(v,w),v+w)\Leftrightarrow B_2(u,w,v)=\phi (w,u)+\phi (w+u,u)+\phi
(v+w,u)+\phi (v+w+u,u)=0$.
\item $L(\phi)$ is a RC-loop $\Leftrightarrow (yz\ast x)xz=y(z\ast xx)
~\forall~x=(a,u),y=(b,v),z=(c,w)\in L(\phi)$. $[(b,v)(c,w)\ast
(a,u)](a,u)=(b,v)[(c,w)\ast (a,u)(a,u)]\Leftrightarrow (b+c+\phi
(v,w)+\phi (v+w,u)+\phi (v+w+u,u),v+w)=(b+c+\phi (u,u)+\phi
(v,w),v+w)\Leftrightarrow B_3(u,v,w)=\phi (u,u)+\phi (v+w,u)+\phi
(v+w+u,u)=0$.
\end{enumerate}

\begin{mycor}\label{1:9}
Let $L(\phi)$ be a code loop on a binary linear code $C$ that is
doubly even with factor set $\phi$. $L(\phi)$ is an RC-loop if and
only if $\phi$ is a lwl factor set.
\end{mycor}
{\bf Proof}\\
Following Theorem~\ref{1:8}, $B_1(u,w)=0\Leftrightarrow \phi
(u+w,u)=\phi (u,u)+\phi (w,u)~\forall~u,w\in C$.

\begin{mycor}\label{1:10}
Let $L(\phi)$ be a code loop on a binary linear code $C$ that is
doubly even with factor set $\phi$. If $L(\phi)$ is an RC-loop ;
\begin{enumerate}
\item $\phi (u+v,u)\equiv [\phi (u,u)+\frac{u\cdot
v}{2}+\phi(u,v)]\bmod{2}$.
\item The difference between $\phi$ been rwl and lwl is $\frac{u\cdot
v}{2}\bmod{2}~\forall~u,v\in C$.
\end{enumerate}
\end{mycor}
{\bf Proof}\\
These follows from Lemma~\ref{1:2} and Corollary~\ref{1:3}.

\begin{mycor}\label{1:11}
Let $L(\phi)$ be a code loop on a binary linear code $C$ that is
doubly even with factor set $\phi$. $L(\phi)$ is a C-loop if and
only if any of the following equivalent conditions is true for all
$u,v,w\in C$.
\end{mycor}
{\bf Proof}\\
\begin{enumerate}
\item $A_1(u,v)=0$ and $B_1(u,v)=0$.
\item $A_2(u,v,w)=0$ and $B_2(u,v,w)=0$.
\item $A_3(u,v,w)=0$ and $B_3(u,v,w)=0$.
\end{enumerate}
{\bf Proof}\\
$L(\phi )$ is a C-loop if and only if it is both a LC-loop and an
RC-loop. The rest follow from Theorem~\ref{1:5} and
Theorem~\ref{1:8}.

\begin{mycor}\label{1:12}
Let $L(\phi)$ be a code loop on a binary linear code $C$ that is
doubly even with factor set $\phi$. $L(\phi)$ is a C-loop if and
only if $\phi$ is a wl factor set on $C$.
\end{mycor}
{\bf Proof}\\
$L(\phi )$ is a C-loop if and only if it is both a LC-loop and an
RC-loop. The rest follow from Corollary~\ref{1:6} and
Corollary~\ref{1:9}.

\begin{mycor}\label{1:13}
Let $L(\phi)$ be a code loop on a binary linear code $C$ that is
doubly even with factor set $\phi$. If $L(\phi)$ is a C-loop ;
\begin{enumerate}
\item $\phi (u,u+v)+\phi (u+v,u)\equiv \frac{u\cdot v}{2}\bmod{2}$.
\item The difference between $\phi$ been rwl and lwl is $\frac{u\cdot
v}{2}\bmod{2}~\forall~u,v\in C$.
\item $\phi (u,u+v)=\phi (u+v,v)\Leftrightarrow \phi (u,v)=\phi
(v,u)~\forall~u,v\in C$.
\end{enumerate}
\end{mycor}
{\bf Proof}\\
Since a C-loop is both an RC-loop and an LC-loop, (1) and (2) follow
from Corollary~\ref{1:10} and Corollary~\ref{1:7}. (3) follows from
Corollary~\ref{1:12} and Lemma~\ref{1:4.2}.

\begin{myth}\label{1:14}
Let $L(\phi)$ be a code loop on a binary linear code $C$ that is
doubly even with factor set $\phi$. $L(\phi)$ is a C-loop if and
only if $D(u,v,w)=\phi (v,u)+\phi (u,w)+\phi (v+u,v)+\phi
(u,u+w)=0$.
\end{myth}
{\bf Proof}\\
$L(\phi)$ is a C-loop $\Leftrightarrow y(x\ast xz)=(yx\ast
x)z~\forall~x=(a,u),y=(b,v),z=(c,w)\in L(\phi)$. $(b,v)[(a,u)\ast
(a,u)(c,w)]=[(b,v)(a,u)\ast (a,u)](c,w)\Leftrightarrow (b+c+\phi
(u,w)+\phi (u,u+w)+\phi (v,w),v+w)=(b+c+\phi (v,u)+\phi (v+u,v)+\phi
(v,w),v+w)\Leftrightarrow D(u,v,w)=\phi (v,u)+\phi (u,w)+\phi
(v+u,v)+\phi (u,u+w)=0$.

\section{Code Loops and Conjugacy Close Loops}
\begin{mylem}\label{1:14.1}
Let $L(\phi)$ be a code loop on a binary linear code $C$ that is
doubly even with factor set $\phi$. $L(\phi)$ is a LCC-loop if and
only if $\phi$ is a rwl factor set.
\end{mylem}
{\bf Proof}\\
A code loop is a Moufang loop, hence an LB-loop. A LB-loop is an
LC-loop if and only if it is an LCC-loop. Thus by
Corollary~\ref{1:6}, the result follows.

\begin{mylem}\label{1:15}
Let $L(\phi)$ be a code loop on a binary linear code $C$ that is
doubly even with factor set $\phi$. $L(\phi)$ is a RCC-loop if and
only if $\phi$ is a lwl factor set.
\end{mylem}
{\bf Proof}\\
A code loop is a Moufang loop, hence an RB-loop. A RB-loop is an
RC-loop if and only if it is an RCC-loop. Thus by
Corollary~\ref{1:9}, the result follows.

\begin{mylem}\label{1:16}
Let $L(\phi)$ be a code loop on a binary linear code $C$ that is
doubly even with factor set $\phi$. $L(\phi)$ is a CC-loop if and
only if $\phi$ is a wl factor set.
\end{mylem}
{\bf Proof}\\
A code loop is a Moufang loop which is possible if and only if it is
a C-loop. Hence, the result follows from Lemma~\ref{1:14.1} and
Lemma~\ref{1:15}.

\section{Code Loops and Extra Loops}
\begin{mylem}\label{1:17}
Let $L(\phi)$ be a code loop on a binary linear code $C$ that is
doubly even with factor set $\phi$. $L(\phi)$ is an extra loop if
and only if $\phi$ is a wl factor set.
\end{mylem}
{\bf Proof}\\
A code loop is a Moufang loop. A Moufang loop is an extra loop if
and only if it is a CC-loop. Hence, the claim follows from
Lemma~\ref{1:16}.

\begin{mycor}\label{1:18}
Let $L(\phi)$ be a code loop on a binary linear code $C$ that is
doubly even with factor set $\phi$. $L(\phi)$ is an extra loop if
and only if $L(\phi )$ is nuclear square.
\end{mycor}
{\bf Proof}\\
A code loop is a Moufang loop. A Moufang loop is an extra loop if
and only if it is nuclear square. Hence, the proof.

\begin{myth}\label{1:20}
Let $L(\phi)$ be a code loop on a binary linear code $C$ that is
doubly even with factor set $\phi$. $L(\phi)$ is an extra loop if
and only if any of the following equivalent conditions is true.
\begin{enumerate}
\item $E_1(u,v,w)=\phi (u,v)+\phi (w,u)+\phi (u+v,w)+\phi
(v+w,u)+\phi (v,w+u)+\phi (u,v+w)=0$.
\item $E_2(u,v,w)=\phi (u,u)+\phi (u,v)+\phi (v,u)+\phi (u,w)+\phi
(v+u,w)+\phi (u,v+w)+\phi (u+v,u+w)=0$.
\item $E_3(u,v,w)=\phi (u,u)+\phi (v,u)+\phi (u,w)+\phi (w,u)+\phi
(v,u+w)+\phi (v+w,u)+\phi (v+u,w+u)=0$.
\end{enumerate}
\end{myth}
{\bf Proof}\\
\begin{enumerate}
\item $L(\phi)$ is a extra loop $\Leftrightarrow (xy\ast z)x=x(y\ast zx)~\forall~x=(a,u),y=(b,v),z=(c,w)\in L(\phi)$. $[(a,u)(b,v)\ast (c,w)](a,u)=(a,u)[(b,v)\ast (c,w)(a,u)]\Leftrightarrow (b+c+\phi
(u,v)+\phi (u+v,w)+\phi (u+v+w,u),v+w)=(b+c+\phi (w,u)+\phi
(v,w+u)+\phi (u,v+w+u),v+w)\Leftrightarrow E_1(u,v,w)=\phi
(u,v)+\phi (w,u)+\phi (u+v,w)+\phi (v+w,u)+\phi (v,w+u)+\phi
(u,v+w)=0$.
\item $L(\phi)$ is a extra loop $\Leftrightarrow xy\ast xz=x(yx\ast z )~\forall~x=(a,u),y=(b,v),z=(c,w)\in L(\phi)$. $(a,u)(b,v)\ast (a,u)(c,w)=(a,u)[(b,v)(a,u)\ast (c,w)]\Leftrightarrow (b+c+\phi
(u,v)+\phi (u,w)+\phi (u+v,u+w),v+w)=(b+c+\phi (v,u)+\phi
(v+u,w)_\phi (u,v+u+w),v+w)\Leftrightarrow E_2(u,v,w)=\phi
(u,u)+\phi (u,v)+\phi (v,u)+\phi (u,w)+\phi (v+u,w)+\phi
(u,v+w)+\phi (u+v,u+w)=0$.
\item $L(\phi)$ is a extra loop $\Leftrightarrow yx\ast zx=(y\ast xz)x~\forall~x=(a,u),y=(b,v),z=(c,w)\in L(\phi)$. $(b,v)(a,u)\ast (c,w)(a,u)=[(b,v)\ast (a,u)(c,w)](a,u)\Leftrightarrow (b+c+\phi
(v,u)+\phi (w,u)+\phi (v+u,w+u),v+w)=(b+c+\phi (u,w)+\phi
(v,u+w)+\phi (v+u+w,u),v+w)\Leftrightarrow E_3(u,v,w)=\phi
(u,u)+\phi (v,u)+\phi (u,w)+\phi (w,u)+\phi (v,u+w)+\phi
(v+w,u)+\phi (v+u,w+u)=0$.
\end{enumerate}

\section{Burn Loops and Code Loops}
\begin{mylem}\label{1:23}
Let $L(\phi)$ be a code loop on a binary linear code $C$ that is
doubly even with factor set $\phi$. $L(\phi)$ is left Burn loop if
and only if $\phi$ is rwl.
\end{mylem}
{\bf Proof}\\
A left Burn loop is a LB-loop that is also an LCC-loop. The result
follows by Lemma~\ref{1:14.1}.

\begin{mylem}\label{1:24}
Let $L(\phi)$ be a code loop on a binary linear code $C$ that is
doubly even with factor set $\phi$. $L(\phi)$ is right Burn loop if
and only if $\phi$ is lwl.
\end{mylem}
{\bf Proof}\\
A right Burn loop is a RB-loop that is also an RCC-loop. The result
follows by Lemma~\ref{1:15}.

\begin{mylem}\label{1:25}
Let $L(\phi)$ be a code loop on a binary linear code $C$ that is
doubly even with factor set $\phi$. $L(\phi)$ is Burn loop if and
only if $\phi$ is wl.
\end{mylem}
{\bf Proof}\\
A Burn loop is a Moufang loop that is also an CC-loop. The result
follows by Lemma~\ref{1:16}.

\begin{mycor}\label{1:26}
Let $L(\phi)$ be a code loop on a binary linear code $C$ that is
doubly even with factor set $\phi$. The following statements about
$L(\phi )$ are equivalent.
\begin{enumerate}
\item $L(\phi )$ is an LC(RC)-loop.
\item $L(\phi )$ is an LCC(RCC)-loop.
\item $L(\phi )$ is a left(right) Burn loop.
\item $\phi$ is rwl(lwl).
\end{enumerate}
\end{mycor}
{\bf Proof}\\
This follows from Corollary~\ref{1:6}, Corollary~\ref{1:9},
Lemma~\ref{1:14.1}, Lemma~\ref{1:15}, Lemma~\ref{1:23} and
Lemma~\ref{1:24}.

\begin{mycor}\label{1:27}
Let $L(\phi)$ be a code loop on a binary linear code $C$ that is
doubly even with factor set $\phi$. The following statements about
$L(\phi )$ are equivalent.
\begin{enumerate}
\item $L(\phi )$ is a C-loop.
\item $L(\phi )$ is a CC-loop.
\item $L(\phi )$ is a Burn loop.
\item  $L(\phi )$ is an extra loop.
\item $\phi$ is wl.
\end{enumerate}
\end{mycor}
{\bf Proof}\\
This follows from Corollary~\ref{1:12}, Lemma~\ref{1:16}, and
Lemma~\ref{1:25} or just by Corollary~\ref{1:26}.

\begin{myth}\label{1:28}
Let $L(\phi)$ be a code loop on a binary linear code $C$ that is
doubly even with factor set $\phi$. The following statements about
$L(\phi )$ are true.
\begin{enumerate}
\item $L(\phi )$ is a C-loop.
\item $L(\phi )$ is a CC-loop.
\item $L(\phi )$ is a Burn loop.
\item  $L(\phi )$ is an extra loop.
\end{enumerate}
\end{myth}
{\bf Proof}\\
$L(\phi )$ is a Moufang loop, so it is alternative which is true if
and only if it is left alternative(LA) and right alternative(RA). It
is easy to show that $L(\phi )$ is LA if and only if $\phi$ is rwl
and $L(\phi )$ is RA if and only if $\phi$ is lwl. So, $L(\phi )$ is
alternative if and only if $\phi$ is wl. So following
Corollary~\ref{1:27}, the claims are true.

\begin{mycor}\label{1:29}
Let $L(\phi)$ be a code loop on a binary linear code $C$ that is
doubly even with factor set $\phi$. The following are true.
\begin{enumerate}
\item  $\phi (u,v)+\phi (w,u)+\phi (u+v,w)+\phi (v,w+u)\equiv
\frac{u\cdot (v+w)}{2}\bmod{2}$.
\item $\phi (u+v,u+w)\equiv [\frac{||u||}{4}+\frac{u\cdot v}{2}+\phi
(u,w)+\phi (v+u,w)+\phi (u,v+w)]\bmod{2}$.
\item $\phi (u+v,u+w)\equiv [\frac{||u||}{4}+\frac{u\cdot w}{2}+\phi
(v,u)+\phi (v,u+w)+\phi (v+w,u)]\bmod{2}$.
\item $\phi (v,u)+\phi (w,u)+\phi (u+v,w)+\phi (v,w+u)\equiv
\frac{u\cdot v+u\cdot (v+w)}{2}\bmod{2}$.
\item $\phi (w,u)+\phi (v,w)+\phi (u,v+w)+\phi (v,w+u)\equiv
[c(u,v,w)+\frac{u\cdot (v+w)}{2}]\bmod{2}$.
\item $\phi (u+v,u)=\phi (u,u)+\phi (v,u)$ and $\phi (u,u+v)=\phi
(u,u)+\phi (u,v)$.
\end{enumerate}
\end{mycor}
{\bf Proof}\\
All these are achieved by Theorem~\ref{1:20}, equations
(\ref{eq:2}), (\ref{eq:3}) and (\ref{eq:4}). We use
Theorem~\ref{1:20} because by Theorem~\ref{1:28}, $L(\phi)$ is an
extra loop.

\begin{myrem}\label{1:30}
The congruence equations in Corollary~\ref{1:29} generalise the
congruence equations (\ref{eq:2}), (\ref{eq:3}) and (\ref{eq:4}).
\end{myrem}

\end{document}